\documentclass[11pt]{article}

\usepackage[letterpaper,width=135mm,height=203mm]{geometry}
\renewcommand{\baselinestretch}{1.25}
\usepackage{graphicx}
\usepackage[normalem]{ulem}
\usepackage{xcolor}

\newtheorem{theorem}{Theorem}
\newtheorem{lemma}[theorem]{Lemma}
\newtheorem{claim}{Claim}

\usepackage{latexsym}
\newenvironment{Proof}{\proofing}{\QED}
\newcommand{\QED}{\hspace{8mm}\mbox{$\Box$}\smallskip}
\newcommand{\proofing}{\textsc{Proof.}}

\newcommand{\curlyC}{\mathcal{C}}
\newcommand{\curlyF}{\mathcal{F}}
\newcommand{\curlyO}{\mathcal{O}}

\newcommand{\dd}[1]{\textbf{\textit{#1}}}

\newcommand{\indiso}{\iota^i}

\makeatletter
\long\def\@caption#1[#2]#3{\begingroup \@parboxrestore 
\if@minipage \@setminipage \fi \normalsize \sffamily \@makecaption {\csname fnum@#1\endcsname }{\ignorespaces #3}\par \endgroup}
\makeatother

\begin{document}


\begin{center} 
\textbf{\Large Bounds on Independent Isolation in Graphs} \\[2mm]
Geoffrey Boyer and Wayne Goddard \\
School of Mathematical and Statistical Sciences, Clemson University
\end{center}

\begin{abstract}
An isolating set of a graph is a set of vertices $S$ such that, if $S$ and its neighborhood is removed,
only isolated vertices remain; and the isolation number is the minimum size of such a set. 
It is known that 
for every connected graph apart from $K_2$ and $C_5$, the isolation number is at most one-third the order and indeed
such a graph has three disjoint isolating sets. In this paper we consider isolating sets where $S$ is 
required to be an independent set and call the minimum size thereof the independent isolation number. 
While for general graphs of order $n$ the independent isolation number 
can be arbitrarily close to $n/2$, we show that in bipartite graphs
the vertex set
can be partitioned into three
disjoint independent isolating sets, whence the independent isolation number is at most $n/3$;
while for $3$-colorable graphs the maximum value of the  independent isolation number is $(n+1)/3$.
We also provide a bound for $k$-colorable graphs.
\end{abstract}

Keywords: graph, isolating set, independent, coloring

\section{Introduction}

Given a family $\curlyF$ of graphs, Caro and Hansberg~\cite{CH} defined an 
\dd{$\curlyF$-isolating set} of graph $G$ as a set $S$ of vertices such that 
$G-N[S]$ contains no copy of any graph in~$\curlyF$, where $N[S]$ denotes the closed neighborhood of $S$ (meaning
$S$ and its neighbors). 
The case where $\curlyF = \{ K_1 \}$ coincides with dominating set.
The case where $\curlyF = \{ K_2 \}$ is now simply called an \dd{isolating set}. 
In~\cite{BGfirst} it was observed that an isolating set
coincides with the concept of \dd{vertex-edge dominating set}, as introduced by Peters~\cite{Peters}. 
Other choices of $\curlyF$ that have been studied include cliques~\cite{BFK,FK} and $K_{1,2}$~\cite{ZW}.

In this paper we study a restriction on the set $S$. Specifically, 
an \dd{independent isolating set} is one where $S$ is itself an independent set.
Further, the \dd{independent isolation number} $\indiso(G)$ of a graph $G$ is 
the minimum size of an independent isolating set. 
This parameter has been investigated as the 
independent vertex-edge domination number, initially by Peters~\cite{Peters} and
Lewis et al.~\cite{LHHF}. More recently,
Boutrig et al.~\cite{BCHH} showed for example that for
claw-free graphs the ordinary and independent numbers are equal, and that for trees the independent vertex-edge domination
number is at most the (ordinary) domination number. The latter bound was later improved by Boutrig and Chellali~\cite{BC} who
showed that for bipartite graphs the independent vertex-edge domination number is at most half the total domination number.
There has also been work in the isolating set formulation; for example,
Favaron and Kaemawichanurat~\cite{FK} examined the relationship between $\curlyF$-isolation and independent $\curlyF$-isolation
especially in claw-free graphs. There has also been work on \dd{total isolation}, where instead the isolating set $S$ is
required to induce a graph without isolates; see for example~\cite{BC,BGH}.

In this paper we focus on bounds and the related question of disjoint independent isolating sets.
In Section~\ref{s:general} we show that there are graphs where the independent isolation number is asymptotically $n/2$, 
where $n$ is the order.
Further, every graph has two disjoint independent isolating sets, but it is NP-complete 
to determine if a graph has three.
In Section~\ref{s:bipartite} we show that a connected bipartite graph of order at least $3$ has three disjoint independent isolating sets.
From this it follows that the independent isolation number of such a graph is at most~$n/3$, which is best possible,
since there are trees whose ordinary isolation number is~$n/3$, as noted in~\cite{CH,KVK}. 
In Section~\ref{s:tripartite} we show that the maximum value for $3$-colorable graphs is $(n+1)/3$.
In Section~\ref{s:colorable} we provide an upper bound for $k$-colorable graphs for larger $k$.
We conclude with some thoughts on future directions.

\section{General Graphs} \label{s:general}

\begin{figure}[h]
\centerline{\includegraphics{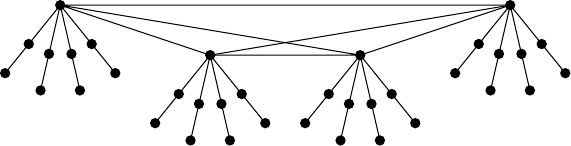}}
\caption{The graph $M_4$}
\label{f:m4}
\end{figure}

A natural question that was not resolved in the literature is: 
what is the maximum independent isolation number of a connected graph for given order $n$.
Peters~\cite{Peters} observed that it is at most $n/2$.
Boutrig et al.~\cite{BCHH} showed that $n/2$ is not attainable for $n\ge 3$ and provided an example with independent isolation number $2n/5-1$.
However, we observe here that it is straightforward to generalize their example to a graph with independent isolation number $n/2 - o(n)$. 
Start with a clique $K_r$, add $r$ feet at each vertex, and then subdivide each leaf-edge once.
The resultant graph $M_r$ has independent isolation number $r(r-1)+1$ while the order is $2r^2+r$;
thus $\indiso(M_r)$ is $n/2 - O( \sqrt{n} )$. Figure~\ref{f:m4} shows the case $r=4$.

The above example shows that there are graphs that do not have three disjoint independent isolating sets.
However, it is immediate that there are two such sets. (Recall that an independent dominating set is an independent set $S$ where every vertex
not in $S$ has a neighbor in $S$.)

\begin{lemma}\label{l: grundy}
Every graph $G$ has two disjoint independent isolating sets. Indeed it has disjoint sets $X$ and $Y$ such that 
$X$ is an independent dominating set of $G$ and $Y$ is an independent isolating set of $G$.
\end{lemma}
\begin{Proof}
Let $X$ be any independent dominating set of $G$. Then let $Y$ be an independent dominating set of the subgraph $G-X$.
Every edge of $G$ has at least one end in $G-X$, and thus has either an end in $Y$ or at least one neighbor in~$Y$.
That is, $Y$ is an isolating set of $G$.
\end{Proof}

\begin{figure}[h]
\centerline{\includegraphics{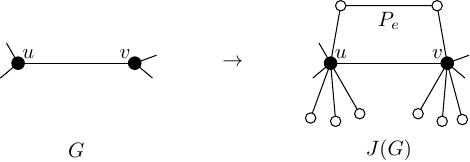}}
\caption{Construction of $J(G)$}
\label{f:jgConstruct}
\end{figure}

We show next that it is NP-complete to determine whether a graph has three disjoint independent isolating sets.
We show this by reducing from chromatic number (which is known to be NP-complete).
For graph $G$, construct graph $J(G)$ from $G$ by:
for each edge $e=uv$ of~$G$, add pair $P_e$ of adjacent vertices with one vertex adjacent to $u$ and the other adjacent to $v$; and
for each vertex $w$ of $G$, add a trio $Q_w$ of end-vertices adjacent only to~$w$. Figure~\ref{f:jgConstruct}
illustrates the construction.

In this paper we often view disjoint independent sets as a partial (proper) coloring of the graph. In such a case,
we say that a vertex is \dd{fully dominated} if its closed neighborhood contains all colors.

\begin{lemma}
For graph $G$, the graph $J(G)$ has three disjoint independent isolating sets if and only if $G$ is $4$-colorable.
\end{lemma}
\begin{Proof}
Suppose $G$ has a proper $4$-coloring. Let $\curlyC$ be a set containing three of these colors. We create
a partial coloring of $J(G)$ using the three colors of $\curlyC$.
First, for vertices of $G$ with color in $\curlyC$, transfer the color to the corresponding vertex of $J(G)$.
Second, for each vertex $w$ of $G$, properly color $Q_w$ such that $Q_w\cup \{w\}$ contains 
all of $\curlyC$. Third, for each edge $e=uv$ of $G$, note that at
least one of $u$ and $v$ has a color in $\curlyC$. So one can properly color both vertices of
$P_e$ such that $P_e \cup \{u,v\}$ contains all of $\curlyC$.

In $J(G)$, each vertex $w$ of $G$ is fully dominated (because of $Q_w$). Further, for each edge $e$ of $G$ 
one vertex of $P_e$ is fully dominated (since both are colored and all of $\curlyC$ is present in $P_e \cup \{u,v\}$). 
So, removing the fully dominated vertices from $J(G)$ leaves
the $Q_w$ and possibly some isolated vertices of the $P_e$; that is, an edgeless graph.
That is, for each color in $\curlyC$ the vertices with that color form  an independent isolating set of $J(G)$.

Conversely, suppose that $J(G)$ has three independent isolating sets, say $(A,B,C)$. 
Then for each edge $e=uv$ of $G$, the set $P_e \cup \{u,v\}$ must 
contain a representative of each of $A$, $B$, and $C$. Thus at least one of $u$ and 
$v$ must be in $A\cup B \cup C$. Furthermore, vertices $u$ and $v$ cannot be from the same set, since they are adjacent in $J(G)$. So we have a proper $4$-coloring of $G$ with color classes
 $A\cap V(G)$, $B\cap V(G)$, $C\cap V(G)$, and $V(G)-(A\cup B \cup C )$ .
\end{Proof}

\section{Bipartite Graphs} \label{s:bipartite}

The \dd{total domination number} of a graph $G$, denoted $\gamma_t(G)$, is the minimum size of a set
$S$ such that every vertex in $G$ has a neighbor in $S$.
Boutrig and Chellali~\cite{BC} established a connection in bipartite graphs between the total domination
number and the independent isolation number:

\begin{theorem} \label{t:bipartiteTDOM}
If a graph $G$ of order $n$ is bipartite, then $\indiso(G) \le \gamma_t  (G) / 2$.
\end{theorem}

For example, since $\gamma_t( G) \le 2n/3$ in any graph with minimum degree at least~$2$ (by \cite{CDH}),
it follows that for such bipartite graphs the independent isolation number is at most $n/3$.
But actually, one can say more.

\begin{theorem}
The vertex set of a connected bipartite graph on at least three vertices can be partitioned into three independent isolating sets.
\end{theorem}
\begin{Proof}
Root the graph $G$ at some vertex $f$, chosen to be an end-vertex of $G$ if there is one. 
This implies that no neighbor of $f$ is an end-vertex.
For each vertex~$w$, let $d_w$ denote the distance from $w$ to $f$.
Color each vertex $w$ by $d_w$ modulo $3$.
Since $G$ is bipartite, each neighbor of $w$ is either at distance $d_w-1$ or
distance $d_w+1$ from $f$;  hence the $3$-coloring is proper. 

It remains to show that each color class
is an isolating set.
If vertex $w$ has neighbors that are both closer to and farther from $f$, then it is immediate that $w$ is fully dominated.
If all neighbors of~$w$ are closer to $f$, then $d_w>1$ since $f$ has no end-vertex neighbor;
thus while $w$ is not dominated by the color of $d_w-2$, all of its neighbors are dominated by that color.
If all neighbors of $w$ are farther, then $w=f$; and since no neighbor of $f$ is an end-vertex, 
each neighbor of $f$ has a neighbor other than $f$ and so is fully dominated.
It follows that deleting the fully dominated vertices destroys all edges, and thus each color class is
an isolating set.
\end{Proof}

As a consequence we get:

\begin{theorem} \label{t:bipartite}
A connected bipartite graph $G$ on $n\ge 3$ vertices has $\indiso(G) \le n/3$.
\end{theorem}

\begin{figure}[h]
\centerline{\includegraphics{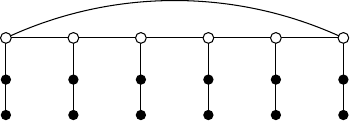}}
\caption{The $P_2$-corona of $C_6$}
\label{f:p2corona}
\end{figure}

The bound of Theorem~\ref{t:bipartite} is best possible because it is best possible for ordinary isolation. 
Caro and Hansberg~\cite{CH} noted that the $P_2$-corona of any graph
has isolation number equal to one-third its order, where the \dd{$P_2$-corona}
means add one end-vertex to each existing vertex, and then subdivide each leaf-edge once.
Figure~\ref{f:p2corona} shows the $P_2$-corona of $C_6$.

For ordinary isolation, the only bipartite graphs with isolation number equal to one-third the order
are the $6$-cycle and the $P_2$-coronas  (see~\cite{BGfirst}). It is unclear if there are more
examples for independent isolation number.

\section{Three-Colorable Graphs} \label{s:tripartite}

In this section we consider $3$-colorable graphs. We show that $(n+1)/3$ is an upper bound
on the independent isolation number and that this value is sharp. Further, such graphs ``almost''
have three disjoint independent isolating sets.

\subsection{Upper Bound}

Consider a proper $3$-coloring $c$ of all the vertices of $G$. We define an edge of $G$ as \dd{bad} if 
the color missing from the edge is also missing from all neighbors of the ends of the edge; otherwise it
is \dd{good}.

\begin{theorem} \label{t:BINGO}
Every connected $3$-colorable graph $G$ has a $3$-coloring where the bad edges, if any, 
have a common vertex. Hence $\indiso( G ) \le (n+1)/3$.
\end{theorem}

\noindent
We prove the theorem in what follows.

\begin{figure}[h]
\centerline{\includegraphics{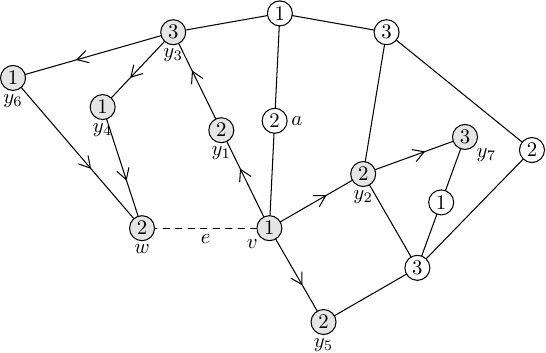}}
\caption{The construction of $D_e$ (shaded) and $F_e$ for $e=vw$}
\label{f:sweep}
\end{figure}

Start with any proper $3$-coloring $c$. Assume there is a bad edge, say $e=vw$.
We build a set $D_e$ of vertices, and a digraph $F_e$ whose vertex set is $D_e$ 
and whose arc set is an orientation of the subgraph of $G-e$ induced by $D_e$.
Assume $c(w)=c(v)+1$ (arithmetic modulo~$3$). Then start with $D_e = \{ v \}$. 
Suppose there is a vertex $y\notin D_e$ such that,
\begin{quote}
(i)~$y$ has at least one neighbor in $D_e$ in $G-e$, \\ 
(ii)~every neighbor of $y$ in $D_e$ has color $c(y)-1$, and\\
(iii)~every neighbor of
$y$ not in $D_e$ has color $c(y)+1$. 
\end{quote}
Then add vertex $y$ to $D_e$, and add to $F_e$ arc(s) to $y$ from 
each neighbor of~$y$ in $D_e$ (but not the arc $wv$).
Repeat this process until there is no such vertex $y$. Figure~\ref{f:sweep} gives an example in which one way to build the $D_e$ is to add vertices in the order $y_1, y_2, \ldots, y_7$, and then add vertex $w$. In contrast, note for instance that vertex $a$ cannot be added since its neighbor outside $D_e$ has color $1$ not color $3$.

\begin{lemma}  \label{l:bingoFD}
At each step, \\
(a) Every arc $xy$ of $F_e$ has $c(y)=c(x)+1$. \\
(b) Every edge $xy$ of $G$ with $x\in D_e$ and $y\notin D_e$ has $c(y)=c(x)+1$.
\end{lemma}
\begin{Proof}
(a) By construction.\\
(b) This is true for $x=v$. For any other vertex $x$, it can only be added to $D_e$ 
if this property is true.
\end{Proof}

By \dd{rotating}  a vertex we mean decrementing its color modulo $3$; that is, changing $1\to 3$, $3\to 2$, or $2\to 1$.
Define the \dd{rotation-sweep out of e} as 
rotating all vertices of $D_e$.

\begin{lemma} \label{l:bingoProper}
After the rotation-sweep out of $e$, the coloring remains proper.
\end{lemma}
\begin{Proof}
An edge in $G$ joining two vertices of
$D_e$ remains proper. An edge of~$G$ joining two vertices not in $D_e$ is unaffected.
By Lemma~\ref{l:bingoFD}, every edge of~$G$ with $x\in D_e$ and $y\notin D_e$ has
$c(y)=c(x)+1$, and so rotating $x$ preserves edge $xy$ being properly colored.
\end{Proof}

\begin{lemma} \label{l:bingoFully}
After the rotation-sweep out of $e$, the following vertices are fully dominated:
(a) any vertex $x\in D_e$ with both an in- and out-arc in $F_e$;
(b) any vertex $y\notin D_e$ that has a neighbor $x\in D_e$ in $G$.
\end{lemma}
\begin{Proof}
(a) This follows since, in $F_e$, in-neighbors have color one less and out-neighbors
have color one more.

(b) 
Assume first that $y\neq w$. Since $y$'s neighbor $x$ was added to $D_e$, it follows that $c(y) = c(x)+1$.
(Note that this is true for all possible $x$.)
Since $y$ was not added to $D_e$, it has a neighbor $z$ not in $D_e$ with the same color as $x$.
So after $x$ is rotated, vertex $y$ has neighbors of both other colors.

Assume second that $y=w$. Suppose $w$ is not 
fully dominated after this rotation-sweep. This means that all neighbors of $w$ were
rotated, and so they were in~$D_e$; thus $w$ should have been added to $D_e$,
a contradiction.
\end{Proof}

\begin{lemma} \label{l:bingoSweep}
After the rotation-sweep out of $e$,\\
(a) if $D_e$ does not contain all of $N_G(v)$, then 
all edges of $G$ incident to a vertex in $N_G[D_e]$ are good; and\\
(b) if $D_e$ does contain all of $N_G(v)$, then
all edges of $G$ incident to a vertex in $N_G[D_e]$ are good except
for edges of the form $vx$ where $x$ is a sink in $F_e$ (meaning no out-neighbor).
\end{lemma}
\begin{Proof}
By definition, neither end of a bad edge is well-dominated.
By  Lemma~\ref{l:bingoFully}, the only possible bad edges with an end in $N_G[D_e]$ 
have both ends being either $v$ or a sink in $F_e$. 

We claim there is no edge of $G$ joining two sinks in $F_e$. 
For, suppose there exists edge $e=xy$ of $G$ where both $x$ and $y$ are sinks in $F_e$.
Say $c(y)=c(x)+1$. Then by construction, when $y$ is added,
vertex $x$ is already in $D_e$. But then by construction, we would have added the arc $xy$ to $F_e$, a contradiction of $x$ being a sink.

So consider an edge joining $v$ to a sink in $F_e$, say $x$. 
Then, if $D_e$ does not contain all of $N_G(v)$, vertex $v$ becomes well-dominated. 
And so part (a) is proved.
But if $D_e$ does contain all of $N_G(v)$, all of $N_G(v)$ remains the same color, 
while all neighbors of $x$ become the same color. Hence the edge $vx$ becomes/stays bad.
\end{Proof}

\begin{lemma} \label{l:bingoRepeat}
If $D_e$ does not contain all vertices in $G$, then repeated rotation-sweeps out of $e$ 
will eventually reach the stage where no edge incident with the original $N_G[D_e]$
is bad.
\end{lemma}
\begin{Proof}
By Lemma~\ref{l:bingoSweep}, we are done immediately unless the initial $D_e$ contains all of $N_G(v)$.
In that case, edge $e$ remains bad. We then do a second rotation-sweep out of $e$. That is, we
construct the set $D'_e$ in the recolored graph and then rotate all vertices in $D'_e$.

We claim that $D'_e$ is strictly contained in $D_e$. For suppose there is a vertex in $D'_e$ but not $D_e$,
and let $y$ be the first such vertex added to $D'_e$. Then $y$ has a neighbor, say $x$, in $D_e \cap D'_e$. So 
originally $c(y) = c(x)+1$, 
but now it is the other way around, which means $y$ cannot be added to $D'_e$ after $x$, a contradiction. 
It follows that $D'_e$ is contained in $D_e$.
Further, by connectivity there must exist a vertex $x\in D_e$ with a neighbor $y\notin D_e$. Then 
since now $c'(x) = c'(y)+1$,  vertex $x$ cannot be added to $D'_e$. The claim follows.

An edge incident with $N_G[D_e]$ but
not $N_G[D'_e]$ remains good, since neither end nor any neighbor of an end changes color.
So again we are done unless $D'_e$ contains all of $N_G(v)$. But if so, we repeat, and 
in a finite number of sweeps, we are done.
\end{Proof}

\textsc{Proof of Theorem~\ref{t:BINGO}.}
Take a $3$-coloring of $G$ with the minimum number of bad edges.
If there is no bad edge then we are done. So
assume there is at least one bad edge, say  $e=vw$ with $c(w)=c(v)+1$. Then do a rotation-sweep
out of~$e$. If $D_e$ does not contain all the vertices of $G$, then by Lemma~\ref{l:bingoRepeat} we
will eventually get rid of the bad edges incident with $v$, and the total number
of bad edges in the graph $G$ will have decreased, a contradiction. So it must be that $D_e$ 
contains all vertices of $G$.

Thus, after the rotation-sweep out of $e$, all bad edges in $G$
are incident with~$v$, and so the first part of the theorem follows. 
Further, the addition of $v$ to the color class of the color missing from $N_G[v]$ produces
an independent isolating set. That is, we have three independent isolating sets whose
total cardinality is $n+1$. The bound follows by averaging.~\QED


\subsection{Optimality}

We show next that the bound of Theorem~\ref{t:BINGO} is sharp.
We will need an operation on a graph $G$:
\begin{quote}
\dd{Operation $\curlyO$.} For vertex $x$, add a vertex $x'$ with the same neighborhood as $x$ (with $x$ and $x'$ nonadjacent), 
and then connect $x'$ to $x$ with a path $xyy'x'$ of length $3$.
\end{quote}

It is immediate that this operation preserves $3$-colorability: give $x'$ the same color as $x$ and give $y$ and $y'$ the other two colors.

\begin{lemma} \label{l:grower}
Operation~$\curlyO$ increases the independent isolation number by exactly one. 
\end{lemma}
\begin{Proof}
Let $G'$ denote the resultant graph.
Any  independent isolating set of $G$ can be extended to one of $G'$ by  
including vertex $y'$. Conversely,  let $J$ be an independent isolating set of $G'$.
If any vertex of $N(x) \cap V(G)$ is in $J$, then $J$ must contain exactly one of $y$ and $y'$, say the latter and the set $J-y'$ independently isolates $G$. 
If no vertex of $N(x)$ is in $J$ and exactly one of $x,x',y',y$ is in $J$, 
then it must be either $y$ or $y'$; say the latter. Then $J-y'$ independently isolates $G$. 
If no vertex of $N(x)$ is in $J$ and two of $x,x',y',y$ are in $J$, then at least one of $x$ or $x'$ must be in $J$, say the former. 
Then $J$ without the other vertex on $P_3$ independently isolates $G$. 
That is, one can always find a vertex of $J$ whose removal leaves an isolating set of $G$.
The claim follows.
\end{Proof}

\begin{figure}[h]
\centerline{\includegraphics{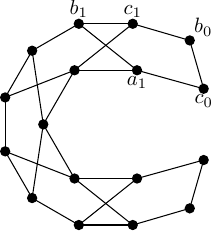}}
\caption{A $3$-colorable graph with $n=17$ and $\indiso =6$}
\label{f:jewel}
\end{figure}

Figure~\ref{f:jewel} shows a specific example.
We build graph $J_m$ as follows. Start with $K_2$ with vertices labeled $b_0c_0$. Then 
apply  Operation~$\curlyO$ to $b_0$, labeling the clone as $a_1$ and the path to $b_0$ having 
interior vertices $b_1$ and $c_1$. Repeat, applying Operation~$\curlyO$ to $b_1$ through $b_{m-2}$.  Figure~\ref{f:jewel} shows the graph $J_6$.
By the above lemma, it holds that $\indiso(J_m) = m$.

\subsection{Cubic Graphs}

\begin{figure}[h]
\centerline{\includegraphics{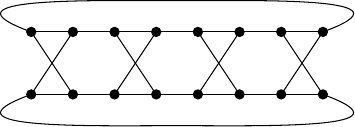}}
\caption{A cubic bipartite graph with $\indiso = n/4$}
\label{f:bipCubMax}
\end{figure}

Archdeacon et al.~\cite{TotalCubic} showed that the total domination number of a cubic graph is at most half its order.
Combined with Theorem~\ref{t:bipartiteTDOM}, it follows that a bipartite cubic
graph has independent isolation number at most one-quarter its order. This is 
best possible, because the extremal cubic graphs for total domination number 
have the desired independent isolation number. Figure~\ref{f:bipCubMax} shows an example.

\begin{figure}[h]
\centerline{\includegraphics{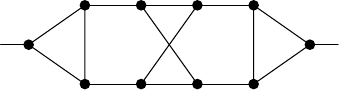}}
\caption{Building block $B$ of cubic graph with $\indiso = \frac{3}{10}n$}
\label{f:cubMax}
\end{figure}

The graph $K_4$ has independent isolation number $1$. Other connected cubic graphs $G$ are $3$-colorable,
and so by Theorem~\ref{t:BINGO} it holds that $\indiso(G) \le (n+1)/3$ where $n$ is the order. 
However, we do not know of an example achieving this bound. For $n>6$ 
the largest value we know of is $\frac{3}{10} n$. Take $r\ge 2$
copies of the $10$-vertex graph $B$ shown in Figure~\ref{f:cubMax}, and join up the copies to form a connected
cubic graph $B_r$. It can readily be argued that an independent isolating set contains 
at least three vertices from every copy of $B$, and that $\indiso(B_r) = 3r$.

\subsection{Maximal Outerplanar Graphs}

One special family of $3$-colorable graphs are the maximal outerplanar graphs.
Caro and Hansberg~\cite{CH} showed that the isolation number of a maximal outerplanar graph is at most $n/4$
and this is sharp. In fact, more is true. (Possibly implicit elsewhere in the literature.) 
Specifically:

\begin{lemma}
Every maximal outerplanar graph of order at least $4$ has 
a partition of the vertex set into four independent isolating sets. 
\end{lemma}

This is an immediate
corollary of  the coloring result of Tokunaga~\cite{Tokunaga}: a maximal outerplanar graph has a proper
coloring with four colors such that every vertex of degree three or more has all
colors in its closed neighborhood. This means that when one removes
all vertices of a particular color along with their neighbors, what remains is a subset of
the vertices of degree $2$; and in a maximal outerplanar graph of order at least $4$ the vertices of degree $2$ form an independent set.


\section{A Bound for $k$-Colorable Graphs} \label{s:colorable}

Finally, we show that for each integer $k\ge 4$ there exists a constant $c_k<\frac{1}{2}$ such that 
for every connected $k$-colorable graph $G$ of order $n$ other than $K_2$ it holds that $\indiso(G) \le c_k n$.

\begin{theorem} \label{t:colorable}
If connected graph $G$  with order $n\ge 3$ is $k$-colorable for $k\ge 4$, then $\indiso( G) \le (k+2)n/(2k+6)$.
\end{theorem}
\begin{Proof}
Consider a Grundy coloring $(A_1, \ldots, A_k)$ of $G$; that is, every vertex of~$A_m$ has a neighbor in each of $A_1$ through $A_{m-1}$.
Let $H$ be the subgraph of~$G$ induced by $A_1 \cup A_2$, and 
let $X$ denote the subgraph of $H$ consisting of the components of order at least $3$.
Let $a_i = |A_i |$, $h=a_1+a_2$, and $x=|V(X)|$.

\begin{claim}
(i) $\indiso \le h/2$. \\
(ii) $\indiso \le n/2 + (n-2h+x)/(2k-4)$.  \\
(iii) $\indiso \le n/2 -x/6$.
\end{claim}
\begin{Proof}
(i) As seen in the proof of Lemma~\ref{l: grundy}, $A_1$ and $A_2$ are independent isolating sets of $G$. 

(ii) 
For each value of $m=3$ through $k$, build an independent isolating set $J_m$ of $G$ starting with $A_m$
and then adding a minimum independent isolating set of each nontrivial component of $G-N[A_m]$,
which we know contains at most half its non-isolated vertices.

Consider a component $F$ of $H-X$; it has order $1$ or $2$. By 
connectivity there is an edge from $F$ to $A_3 \cup \ldots \cup A_k$. Let $m^*\ge 3$ be the smallest value such 
that there is an edge from $F$ to $A_{m^*}$. The subgraph
$G-N[A_{m^*}]$ contains at most one vertex of $F$ and no vertex of $A_{m^*+1}$ through $A_k$ (because of the Grundy coloring);
hence if one vertex of $F$ remains it is isolated. 
It follows that the sum of the numbers of non-isolated vertices in $G-N[A_3]$ through $G-N[A_k]$  is at most
$ \left( \sum_{m=3}^{k}  (n - a_m) \right)  - (h-x) = (k-2)n - (n-h) - (h-x) = (k-2)n - n+x$. Hence 
the total size of the independent isolating sets $J_3 $ through $J_k$ is at most $(n-h) + ((k-2) n - n+x) )/2 = (k-2)n/2 + (n-2h+x)/2$.
The result follows by averaging. (Aside: this bound can clearly be improved but that doesn't 
help the overall bound.)

(iii) Build an independent isolating set of $G$ starting with a minimum 
independent isolating set of $X$. Since $X$ is bipartite, 
this uses at most $x/3$ vertices, by Theorem~\ref{t:bipartite}. Then, as before, we need
at most half the remaining vertices. Thus $\indiso \le x/3 + (n-x)/2 = n/2 - x/6$.
\end{Proof}


The above claim yields a linear program; namely, maximize $f$, say, subject to the constraint that 
$f$ is at most each of the three bounds.
This linear program has an optimum of $(k+2)n/(2k+6)$.
This can be checked by noting that $(h,x,f) = \frac{n}{2k+6} ( 2k+4, 6,  k+2 )$
is feasible;  that 
$\frac{n}{k+3} ( 2, k-2, 3 )$ satisfies the dual;
and both have the stated value.
\end{Proof}

It is unclear what the actual maximum independent isolation number is. But the graphs $M_k$ have 
independent isolation number approximately $n( \frac{1}{2} - \frac{1}{2k} )$, which is not too far off the bound
of Theorem~\ref{t:colorable}.



\section{Future Questions}

There remain many questions. One direction would be to determine better bounds for families of graphs;
for example with a floor on the minimum degree. Another direction is to continue the study of independent $\curlyF$-isolation for other choices of $\curlyF$.


\renewcommand{\baselinestretch}{1}
\normalfont

\end{document}